\documentclass{article}
\usepackage{amssymb}
\usepackage{amsfonts}
\usepackage{amsmath}
\usepackage[left=1.5in, right=1.5in, bottom=1.5in,top=1.5in]{geometry}

\input{tcilatex}
\begin{document}

\author{Mohamed Ibrahim$^{1},$ M. S. Eliwa$^{2,}$\thanks{%
Corresponding author: mseliwa@mans.edu.eg} , M. El- Morshedy$^{2}$ \\
\\
$^{1}$Department of Applied Statistics and Insurance, \\
Faculty of Commerce, Damietta University, Damietta, Egypt.\\
$^{2}$Department of Mathematics, Faculty of Science, \\
Mansoura University, Mansoura, Egypt.\\
\\
mohamedibrahim1082@gmail.com\\
mseliwa@mans.edu.eg\\
mah\_elmorshedy@mans.edu.eg\\
}
\title{Bivariate Exponentiated Generalized Linear Exponential Distribution
with Applications in Reliability Analysis}
\date{}
\maketitle

\begin{abstract}
The aim of this paper, is to define a bivariate exponentiated generalized
linear exponential distribution based on Marshall-Olkin shock model.
Statistical and reliability properties of this distribution are discussed.
This includes quantiles, moments, stress-strength reliability, joint
reliability function, joint reversed (hazard) rates functions and joint mean
waiting time function. Moreover, the hazard rate, the availability and the
mean residual lifetime functions for a parallel system, are established. One
data set is analyzed, and it is observed that, the proposed distribution
provides a better fit than Marshall-Olkin bivariate exponential, bivariate
generalized exponential and bivariate generalized linear failure rate
distributions. Simulation studies are presented to estimate both the
relative absolute bias, and the relative mean square error for the
distribution parameters based on complete data.

\end{abstract}

\ \ \ \ \ \ \ \ \ \ \ 

\section{Introduction}

Sarhan et al. (2013) introduced exponentiated generalized linear exponential
distribution (EGLED), which generalized a lot of probability distributions
such as exponential (E), generalized exponential (GE), linear exponential
(LE), generalized linear failure rate (GLFR), generalized linear exponential
(GLE) distributions, among others. Furthermore, the EGLED provides more
flexibility to analyze real data sets such as Leukemia data, drug data,
among others.

In many scientific practical situations, multivariate lifetime data arise
frequently, so it is important to consider different multivariate models
that could be used to model such multivariate lifetime data. Such these
models are interesting in several applications, such as reliability
engineering, industrial engineering and computer systems. So, the aim of
this paper is to introduce a bivariate exponentiated generalized linear
exponential distribution (BEGLED) based on Marshall-Olkin shock model
(1967), whose marginal distributions are EGLED. In the mentioned
applications of the bivariate distribution, could be the lifetimes of two
components, the magnitudes of stress and strength components and drought
intensities. A lot of bivariate distributions based on Marshall-Olkin model
are studied by many authors, see Sarhan and Balakrishnan (2007), Al-
Khedhairi and El-Gohary (2008), Kundu and Gupta (2009), Sarhan et al.
(2011), Kundu and Gupta (2013), Balakrishna and Shiji (2014), El-Gohary et
al. (2016), Rasool and Akbar (2016) and El-Bassiouny et al. (2016).

The random variable $X$ \ is said to have EGLED$(a,b,\alpha ,\theta )$ if
its CDF is \ \ \ 
\begin{equation}
F_{X}(x)=\left( 1-e^{-\eta ^{\alpha }(x)}\right) ^{\theta };\ \ x\geq 0,
\label{1}
\end{equation}%
where $\eta (x)=ax+\frac{b}{2}x^{2}$, the parameters $a,b\geq 0\ $such that
\bigskip $a+b>0$ and\ $\alpha ,\theta >0$. \ \ \ \ \ \ \ \ \ \ \ \ \ \ \ \ \
\ \ \ \ The parameters $a$ and $b$ are scale parameters, while $\alpha $ and 
$\theta $ are shape parameters.

\section{\textbf{The B}EGLED and Its Marginal Functions}

Assume $U_{i}\sim EGLED(\alpha ,a,b,\theta _{i}),\ i=1,2,3\ $are three
independent random variables. Define $X_{k}=max\{U_{k},U_{3}\}$ $;k=1,2$.
So, the bivariate vector $(X_{1},X_{2})$ has the BEGLED with parameters
vector $\Phi =$($\alpha ,a,b,\theta _{1},\theta _{2},\theta _{3}$). The
joint CDF of $(X_{1},X_{2})$ is 
\begin{equation}
F_{X_{1},X_{2}}(x_{1},x_{2})=\left( 1-e^{-\eta ^{\alpha }(x_{1})}\right)
^{\theta _{1}}\left( 1-e^{-\eta ^{\alpha }(x_{2})}\right) ^{\theta
_{2}}\left( 1-e^{-\eta ^{\alpha }(z)}\right) ^{\theta _{3}},\ \ z=\min
(x_{1},x_{2}).\   \label{2}
\end{equation}%
Also, we can get the joint PDF of $(X_{1},X_{2})$ as follows

\begin{equation}
f_{X_{1},X_{2}}(x_{1},x_{2})=\left\{ 
\begin{array}{l}
f_{1}(x_{1},x_{2})\ \ \ \ \ \ \text{if \ }0<x_{1}<x_{2} \\ 
f_{2}(x_{1},x_{2})\ \ \ \ \ \ \text{if \ }0<\text{\ }x_{2}<x_{1} \\ 
f_{3}(x,x)\ \ \ \ \ \ \ \ \ \text{if\ \ }x_{1}=x_{2}=x,%
\end{array}%
\right.  \label{3}
\end{equation}%
where

\begin{equation*}
f_{1}(x_{1},x_{2})=\varphi _{2}\ \left[ \eta (x_{1})\eta (x_{2})\right]
^{\alpha -1}e^{-(\eta ^{\alpha }(x_{1})+\eta ^{\alpha }(x_{2}))}\left(
1-e^{-\eta ^{\alpha }(x_{2})}\right) ^{\theta _{2}-1}\left( 1-e^{-\eta
^{\alpha }(x_{1})}\right) ^{\theta _{1}+\theta _{3}-1},
\end{equation*}

\begin{equation*}
f_{2}(x_{1},x_{2})=\varphi _{1}\ \left[ \eta (x_{1})\eta (x_{2})\right]
^{\alpha -1}e^{-(\eta ^{\alpha }(x_{1})+\eta ^{\alpha }(x_{2}))}\left(
1-e^{-\eta ^{\alpha }(x_{1})}\right) ^{\theta _{1}-1}\left( 1-e^{-\eta
^{\alpha }(x_{2})}\right) ^{\theta _{2}+\theta _{3}-1},
\end{equation*}%
and

\begin{equation*}
f_{3}(x,x)=\alpha \theta _{3}(a+bx)\eta ^{\alpha -1}(x)e^{-\eta ^{\alpha
}(x)}\left( 1-e^{-\eta ^{\alpha }(x)}\right) ^{\theta _{1}+\theta
_{2}+\theta _{3}-1},
\end{equation*}%
where

\begin{equation*}
\varphi _{i}=\alpha ^{2}\theta _{i}(\theta _{3-i}+\theta
_{3})(a+bx_{1})(a+bx_{2})\ \ \ ,\ \ \ \ i=1,2.
\end{equation*}%
On the other hand, the marginal CDFs for the BEGLED can be represented as
follows 
\begin{equation}
F_{X_{i}}(x_{i})=\left( 1-e^{-\eta ^{\alpha }(x_{i})}\right) ^{\theta
_{i}+\theta _{3}},\ \ i=1,2.  \label{4}
\end{equation}%
Therefore, we can get the marginal PDFs for the BEGLED as follows%
\begin{equation}
f(x_{i})=\alpha (\theta _{i}+\theta _{3})(a+bx_{i})\eta ^{\alpha
-1}(x_{i})e^{-\eta ^{\alpha }(x_{i})}\left( 1-e^{-\eta ^{\alpha
}(x_{i})}\right) ^{\theta _{i}+\theta _{3}-1},\ \ i=1,2.  \label{10}
\end{equation}

\section{Statistical Properties}

\subsection{The median of the BEGLED}

Domma (2009) presented the median correlation coefficient $M_{X_{1},X_{2}}\ $%
as a form \ \ \ \ \ 

\begin{equation}
M_{X_{1},X_{2}}=4F_{X_{1},X_{2}}(M_{X_{1}},M_{X_{2}})-1,
\end{equation}%
where $M_{X_{1}}$ and $M_{X_{2}}$ denote the median of $X_{1}$ and $X_{2}$
respectively.

If $X_{1}\sim EGLED(\alpha ,a,b,\theta _{1}+\theta _{3})\ $and $X_{2}\sim
EGLED(\alpha ,a,b,\theta _{2}+\theta _{3})\ $then%
\begin{equation}
M_{X_{i}}=\frac{1}{b}\left( -a+\sqrt{a^{2}+2b\left[ -\ln \left( 1-\left( 
\frac{1}{2}\right) ^{(\frac{1}{\theta _{i}+\theta _{3}})}\right) \right] ^{%
\frac{1}{\alpha }}}\right) ,\ \ \ \ i=1,2.
\end{equation}%
So, the coefficient of median correlation between $X_{1}$ and $X_{2}$ is

\begin{equation}
M_{X_{1},X_{2}}=\left\{ 
\begin{array}{l}
4\left( 1-e^{-\eta ^{\alpha }(M_{X_{2}})}\right) ^{\theta _{2}}\left(
1-e^{-\eta ^{\alpha }(M_{X_{1}})}\right) ^{\theta _{1}+\theta _{3}}-1\ \ \ \ 
\text{if \ }x_{1}<x_{2} \\ 
4\left( 1-e^{-\eta ^{\alpha }(M_{X_{1}})}\right) ^{\theta _{1}}\left(
1-e^{-\eta ^{\alpha }(M_{X_{2}})}\right) ^{\theta _{2}+\theta _{3}}-1\ \ \ \ 
\text{if \ }x_{1}>x_{2}.%
\end{array}%
\right.  \label{5}
\end{equation}%
Equation (\ref{5}) can be used to generate a bivariate data.

\subsection{The mathematical expectation}

We can derive the marginal expectation ($rth$ moment) of $X_{i}\ $ when $%
X_{i}\sim EGLED(\alpha ,a,b,\theta _{i}+\theta _{3})$ such that $i=1,2$ as
follows 
\begin{equation}
E(X_{i}^{r}\ )=\dint\limits_{0}^{\infty }x_{i}^{r}f_{X_{i}}(x_{i})dx_{i},
\end{equation}%
by using Equation (\ref{10}), Maclaurin expansion, binomial expansion and
gamma function, we get

\begin{equation}
E(X_{i}^{r}\ )=\sum_{j=0}^{\infty }\sum_{k=0}^{\alpha -1}\sum_{l=0}^{\infty
}\zeta _{j,k}^{(l)}(a+b\Theta )\Gamma (\Theta ),\ 
\end{equation}%
where

\begin{eqnarray*}
\zeta _{j,k}^{(l)} &=&\alpha \left( \theta _{i}+\theta _{3}\right) \frac{%
(-1)^{j}a^{k}U_{j}^{(l)}}{l!}\left( \frac{b}{2}\right) ^{\alpha -k-1}\binom{%
\alpha -1}{k}\binom{\theta _{i}+\theta _{3}-1}{j}, \\
\Theta &=&r+l+2\alpha -k-1,
\end{eqnarray*}%
and 
\begin{equation*}
U_{j}^{(l)}=\frac{d^{l}}{dx_{i}^{l}}(\exp (-(j+1)\eta ^{\alpha
}(x_{i})+x_{i}))|_{x_{i}=0}.
\end{equation*}

\subsection{The conditional probability density functions}

The conditional probability density function of $X_{i}$ given $X_{j}=x_{j},$ 
$(i,j=1,2,i\neq j)$ is given by%
\begin{equation}
f_{X_{i}\mid X_{j}}(x_{i}\mid x_{j})=\left\{ 
\begin{array}{l}
f_{X_{i}\mid X_{j}}^{(1)}(x_{i}\mid x_{j})\ \ \ \ \text{if }\ x_{i}>x_{j}>0
\\ 
f_{X_{i}\mid X_{j}}^{(2)}(x_{i}\mid x_{j})\ \ \ \ \text{if \ }x_{j}>x_{i}>0
\\ 
f_{X_{i}\mid X_{j}}^{(3)}(x_{i}\mid x_{j})\ \ \ \ \text{if\ \ }x_{i}=x_{j}>0,%
\end{array}%
\right.  \label{11}
\end{equation}%
where%
\begin{eqnarray*}
f_{X_{i}\mid X_{j}}^{(1)}(x_{i} &\mid &x_{j})=\alpha \theta
_{i}(a+bx_{i})\eta ^{\alpha -1}(x_{i})e^{-\eta ^{\alpha }(x_{i})}\left(
1-e^{-\eta ^{\alpha }(x_{i})}\right) ^{\theta _{i}-1}, \\
f_{X_{i}\mid X_{j}}^{(2)}(x_{i} &\mid &x_{j})=\frac{\alpha \theta
_{j}(\theta _{i}+\theta _{3})(a+bx_{i})\eta ^{\alpha -1}(x_{i})e^{-\eta
^{\alpha }(x_{i})}\left( 1-e^{-\eta ^{\alpha }(x_{i})}\right) ^{\theta
_{i}+\theta _{3}-1}}{(\theta _{j}+\theta _{3})(1-e^{-\eta ^{\alpha
}(x_{j})})^{\theta _{3}}},
\end{eqnarray*}%
and%
\begin{equation*}
f_{X_{i}\mid X_{j}}^{(3)}(x_{i}\mid x_{j})\ =\frac{\theta _{3}}{\theta
_{j}+\theta _{3}}\left[ 1-e^{-\eta ^{\alpha }(x_{i})}\right] ^{\theta _{i}}.
\end{equation*}%
Equation (\ref{11}) can be getting by substituting from Equations (\ref{3})
and (\ref{10}) in the following relation%
\begin{equation}
f_{X_{i}\mid X_{j}}(x_{i}\mid x_{j})\ =\frac{f_{X_{i},X_{j}}(x_{i},x_{j})}{%
f_{X_{j}}(x_{j})\ },\text{ \ }(i\neq j=1,2).
\end{equation}

\subsection{The distributions of $T=\max (X_{1},X_{2})$ and $S=\min
(X_{1},X_{2})$}

In the mentioned applications $X_{1}$ and $X_{2}$ could be exchange rates in
two time periods. So, it is important to get the distributions of $T\ $and $%
S $. If the bivariate vector $(X_{1},X_{2})$ has the BEGLED then 
\begin{eqnarray}
F_{T}(t) &=&P(\max (X_{1},X_{2})\leq t)  \notag \\
&=&P(\max (U_{1},U_{3})\leq t,\max (U_{2},U_{3})\leq t)  \notag \\
&=&F_{EGLED}(\alpha ,a,b,\theta _{1}+\theta _{2}+\theta _{3})\text{.}
\end{eqnarray}

Also, we can get the distribution of $S$ as follows

\begin{eqnarray}
F_{S}(t) &=&P(\min (X_{1},X_{2})\leq t)  \notag \\
&=&P(X_{1}<t)+P(X_{2}<t)-P(X_{1}<t,X_{2}<t,)  \notag \\
&=&F_{EGLED}(\alpha ,a,b,\theta _{1}+\theta _{3})+F_{EGLED}(\alpha
,a,b,\theta _{2}+\theta _{3})  \notag \\
&&-F_{EGLED}(\alpha ,a,b,\theta _{1}+\theta _{2}+\theta _{3})\text{.}
\label{6}
\end{eqnarray}

\section{Reliability Properties}

In this section, we present the stress-strength reliability, the joint
reliability function, the joint reversed (hazard) functions and the joint of
mean waiting time function. Also, we present the hazard rate, the
availability and the mean residual lifetime functions for a parallel system
with two components.

\subsection{Stress-strength reliability}

Let $X_{1}$ is a random variable represents stress, and $X_{2}$ is a random
variable represents strength, and the random vector $(X_{1},X_{2})$ has the
BEGLED then, the reliability function $R\ $is 
\begin{eqnarray}
R &=&P[X_{1}<X_{2}]  \notag \\
&=&P(U_{1}<U_{3}<U_{2})+P(U_{3}<U_{1}<U_{2})  \notag \\
&=&\frac{\theta _{2}+\theta _{3}}{\theta _{1}+\theta _{2}+2\theta _{3}}.
\end{eqnarray}

\subsection{The joint reliability function}

Assume ($X_{1},X_{2}$) be two dimensional random variable with CDF $%
F_{X_{1},X_{2}}(x_{1},x_{2})$,$\ $and the marginal functions are $%
F_{X_{1}}(x_{1})\ $and\ $F_{X_{2}}(x_{2})\ $then, the joint reliability
function $R_{X_{1},X_{2}}(x_{1},x_{2})\ $is%
\begin{equation}
R_{X_{1},X_{2}}(x_{1},x_{2})=1-F_{X_{1}}(x_{1})-F_{X_{2}}(x_{2})+F_{X_{1},X_{2}}(x_{1},x_{2}).
\label{99}
\end{equation}%
Assume the random vector $(X_{1},X_{2})$ has the BEGLED then, the joint
reliability function of $(X_{1},X_{2})$ is given by%
\begin{equation}
R_{X_{1},X_{2}}(x_{1},x_{2})=\left\{ 
\begin{array}{l}
R_{1}(x_{1},x_{2})\ \ \ \ \text{if \ }0<x_{1}<x_{2} \\ 
R_{2}(x_{1},x_{2})\ \ \ \ \text{if \ }0<x_{2}<x_{1} \\ 
R_{3}(x,x)\ \ \ \ \ \ \ \ \text{if\ \ }x_{1}=x_{2}=x,%
\end{array}%
\right.  \label{100}
\end{equation}%
where

\begin{eqnarray*}
R_{1}(x_{1},x_{2}) &=&1-\left( 1-e^{-\eta ^{\alpha }(x_{1})}\right) ^{\theta
_{1}+\theta _{3}}-\left( 1-e^{-\eta ^{\alpha }(x_{2})}\right) ^{\theta
_{2}+\theta _{3}} \\
&&+\left( 1-e^{-\eta ^{\alpha }(x_{2})}\right) ^{\theta _{2}}\left(
1-e^{-\eta ^{\alpha }(x_{1})}\right) ^{\theta _{1}+\theta _{3}},
\end{eqnarray*}

\begin{eqnarray*}
R_{2}(x_{1},x_{2}) &=&1-\left( 1-e^{-\eta ^{\alpha }(x_{1})}\right) ^{\theta
_{1}+\theta _{3}}-\left( 1-e^{-\eta ^{\alpha }(x_{2})}\right) ^{\theta
_{2}+\theta _{3}} \\
&&+\left( 1-e^{-\eta ^{\alpha }(x_{1})}\right) ^{\theta _{1}}\left(
1-e^{-\eta ^{\alpha }(x_{2})}\right) ^{\theta _{2}+\theta _{3}},
\end{eqnarray*}

\begin{eqnarray*}
R_{3}(x,x) &=&1-\left( 1-e^{-\eta ^{\alpha }(x)}\right) ^{\theta _{1}+\theta
_{3}}-\left( 1-e^{-\eta ^{\alpha }(x)}\right) ^{\theta _{2}+\theta _{3}} \\
&&+\left( 1-e^{-\eta ^{\alpha }(x)}\right) ^{\theta _{1}+\theta _{2}+\theta
_{3}}.
\end{eqnarray*}

\subsection{The joint reversed (hazard) rate functions}

\subsubsection{\protect\bigskip The joint hazard rate function and its
marginal functions}

Assume ($X_{1},X_{2}$) be two dimensional random variable with PDF $%
f_{X_{1},X_{2}}(x_{1},x_{2})$,$\ $and reliability function $%
R_{X_{1},X_{2}}(x_{1},x_{2})$. Basu (1971) defined the bivariate hazard rate
function as%
\begin{equation}
h(x_{1},x_{2})=\frac{f_{X_{1},X_{2}}(x_{1},x_{2})}{%
R_{X_{1},X_{2}}(x_{1},x_{2})}.  \label{101}
\end{equation}%
So, the bivariate hazard rate function for the random vector $(X_{1},X_{2})$
which has the BEGLED is \  
\begin{equation}
h_{X_{1},X_{2}}(x_{1},x_{2})=\left\{ 
\begin{array}{l}
h_{1}(x_{1},x_{2})\ \ \ \ \ \ \text{if \ }0<x_{1}<x_{2} \\ 
h_{2}(x_{1},x_{2})\ \ \ \ \ \ \text{if \ }0<x_{2}<x_{1} \\ 
h_{3}(x,x)\ \ \ \ \ \ \ \ \text{if\ \ }x_{1}=x_{2}=x,%
\end{array}%
\right.   \label{102}
\end{equation}%
where

\begin{equation*}
h_{1}(x_{1},x_{2})=\frac{\varphi _{2}\ \left[ \eta (x_{1})\eta (x_{2})\right]
^{\alpha -1}e^{-(\eta ^{\alpha }(x_{1})+\eta ^{\alpha }(x_{2}))}(\Psi
(x_{2}))^{\theta _{2}-1}\left( \Psi (x_{1})\right) ^{\theta _{1}+\theta
_{3}-1}}{1-\left( \Psi (x_{1})\right) ^{\theta _{1}+\theta _{3}}-\left( \Psi
(x_{2})\right) ^{\theta _{2}+\theta _{3}}+\left( \Psi (x_{2})\right)
^{\theta _{2}}\left( \Psi (x_{1})\right) ^{\theta _{1}+\theta _{3}}},
\end{equation*}

\begin{equation*}
h_{2}(x_{1},x_{2})=\frac{\varphi _{1}\ \left[ \eta (x_{1})\eta (x_{2})\right]
^{\alpha -1}e^{-(\eta ^{\alpha }(x_{1})+\eta ^{\alpha }(x_{2}))}\left( \Psi
(x_{1})\right) ^{\theta _{1}-1}\left( \Psi (x_{2})\right) ^{\theta
_{2}+\theta _{3}-1}}{1-\left( \Psi (x_{1})\right) ^{\theta _{1}+\theta
_{3}}-\left( \Psi (x_{2})\right) ^{\theta _{2}+\theta _{3}}+\left( \Psi
(x_{1})\right) ^{\theta _{1}}\left( \Psi (x_{2})\right) ^{\theta _{2}+\theta
_{3}}},
\end{equation*}%
and 
\begin{equation*}
h_{3}(x,x)=\frac{\alpha \theta _{3}(a+bx)\eta ^{\alpha -1}(x)e^{-\eta
^{\alpha }(x)}\left( \Psi (x)\right) ^{\theta _{1}+\theta _{2}+\theta _{3}-1}%
}{1-\left( \Psi (x)\right) ^{\theta _{1}+\theta _{3}}-\left( \Psi (x)\right)
^{\theta _{2}+\theta _{3}}+\left( \Psi (x)\right) ^{\theta _{1}+\theta
_{2}+\theta _{3}}},
\end{equation*}%
where $\Psi (.)=1-e^{-\eta ^{\alpha }(.)}.$

Also, the marginal hazard rate functions $h_{i}(x_{i}),i=1,2$ of the BEGLED
are

\begin{equation}
h_{i}(x_{i})=\frac{\alpha (\theta _{i}+\theta _{3})(a+bx_{i})\eta ^{\alpha
-1}(x_{i})e^{-\eta ^{\alpha }(x_{i})}\left( \Psi (x_{i})\right) ^{\theta
_{i}+\theta _{3}-1}}{1-\left( \Psi (x_{i})\right) ^{\theta _{i}+\theta _{3}}}%
,\ i=1,2.
\end{equation}

\subsubsection{The joint reversed hazard rate function and its marginal
functions}

Assume ($X_{1},X_{2}$) be two dimensional random variable with CDF $%
F_{X_{1},X_{2}}(x_{1},x_{2})$, the joint reversed hazard rate function is 
\begin{equation}
r(x_{1},x_{2})=\frac{f_{X_{1},X_{2}}(x_{1},x_{2})}{%
F_{X_{1},X_{2}}(x_{1},x_{2})}.  \label{103}
\end{equation}%
So, the joint reversed hazard rate function for the random vector $%
(X_{1},X_{2})$ which has the BEGLED is 
\begin{equation}
r_{X_{1},X_{2}}(x_{1},x_{2})=\left\{ 
\begin{array}{l}
r_{1}(x_{1},x_{2})\ \ \ \ \ \ \text{if \ }0<x_{1}<x_{2} \\ 
r_{2}(x_{1},x_{2})\ \ \ \ \ \ \text{if \ }0<x_{2}<x_{1} \\ 
r_{3}(x,x)\ \ \ \ \ \ \ \ \text{if\ \ }x_{1}=x_{2}=x,%
\end{array}%
\right.  \label{104}
\end{equation}%
where

\begin{equation*}
r_{1}(x_{1},x_{2})=\frac{\varphi _{2}\ \left[ \eta (x_{1})\eta (x_{2})\right]
^{\alpha -1}e^{-(\eta ^{\alpha }(x_{1})+\eta ^{\alpha }(x_{2}))}}{\left(
1-e^{-\eta ^{\alpha }(x_{2})}\right) \left( 1-e^{-\eta ^{\alpha
}(x_{1})}\right) },
\end{equation*}

\begin{equation*}
r_{2}(x_{1},x_{2})=\frac{\varphi _{1}\ \left[ \eta (x_{1})\eta (x_{2})\right]
^{\alpha -1}e^{-(\eta ^{\alpha }(x_{1})+\eta ^{\alpha }(x_{2}))}}{\left(
1-e^{-\eta ^{\alpha }(x_{2})}\right) \left( 1-e^{-\eta ^{\alpha
}(x_{1})}\right) },
\end{equation*}%
and 
\begin{equation*}
r_{3}(x,x)=\frac{\alpha \theta _{3}(a+bx)\eta ^{\alpha -1}(x)e^{-\eta
^{\alpha }(x)}}{\left( 1-e^{-\eta ^{\alpha }(x)}\right) }.
\end{equation*}%
Also, the marginal reversed hazard rate functions $r_{i}(x_{i}),i=1,2$ to
the BEGLED are 
\begin{equation}
r_{i}(x_{i})=\frac{\alpha (\theta _{i}+\theta _{3})(a+bx_{i})\eta ^{\alpha
-1}(x_{i})e^{-\eta ^{\alpha }(x_{i})}}{\left( 1-e^{-\eta ^{\alpha
}(x_{i})}\right) },\ i=1,2.  \label{105}
\end{equation}

\subsection{The joint mean waiting time and its marginal functions}

The waiting time is closely related to important random variable reversed
hazard rate function, which the failure occurs in the interval [0, t]. The
observations of waiting times can be used for prediction the distribution
function. So, one of the most important applications of the waiting time is
to describe different maintenance strategies to any system. The joint mean
waiting time function $M_{w}(t_{1},t_{2})$ is defined as follows 
\begin{equation}
M_{w}(t_{1},t_{2})=\frac{1}{F(t_{1},t_{2})}\int_{0}^{t_{1}}%
\int_{0}^{t_{2}}F(x_{1},x_{2})\ dx_{2}dx_{1}.  \label{106}
\end{equation}%
Assume the random vector $(X_{1},X_{2})$ has the BEGLED. Using Maclaurin and
binomial expansions then, the joint mean waiting time function $%
M_{w}(t_{1},t_{2})$ is 
\begin{equation}
M_{w}(t_{1},t_{2})=\left\{ 
\begin{array}{l}
M_{w_{1}}(t_{1},t_{2})\ \ \ \ \ \ \ \ \ \text{if \ }t_{1}>t_{2}>0 \\ 
M_{w_{2}}(t_{1},t_{2})\ \ \ \ \ \ \ \ \ \text{if \ }0<t_{1}<t_{2} \\ 
M_{w_{3}}(t,t)\ \ \ \ \ \ \ \text{if\ \ }t_{1}=t_{2}=t,%
\end{array}%
\right.  \label{107}
\end{equation}%
where 
\begin{equation*}
M_{w_{i}}(t_{1},t_{2})=\frac{1}{F(t_{1},t_{2})}\dsum\limits_{j,k=0}^{\infty }%
\frac{g_{j}^{(k)}\psi _{j}^{(k)}}{\left( (k+1)!\right) ^{2}}\binom{\theta
_{i}}{j}\binom{\theta _{3-i}+\theta _{3}}{j}\left( t_{1}t_{2}\right)
^{k+1};\ \ i=1,2,
\end{equation*}%
\begin{equation*}
M_{w_{3}}(t,t)=\frac{1}{F(t,t)}\dsum\limits_{j,k=0}^{\infty }\frac{%
(-1)^{j}Q_{j}^{(k)}}{(k+1)!}\binom{\theta _{1}+\theta _{2}+\theta _{3}}{j}%
t^{k+1},
\end{equation*}%
and 
\begin{equation*}
g_{j}^{(k)}=\frac{d^{k}}{dx_{1}^{k}}(e^{-j\eta ^{\alpha
}(x_{1})})|_{x_{1}=0}\ ,\ \ \psi _{j}^{(k)}=\frac{d^{k}}{dx_{2}^{k}}%
(e^{-j\eta ^{\alpha }(x_{2})})|_{x_{2}=0},\text{\ }\ \ Q_{j}^{(k)}=\frac{%
d^{k}}{dx^{k}}(e^{-j\eta ^{\alpha }(x)})|_{x=0}.
\end{equation*}%
Also, the marginal mean waiting time functions $m_{w_{i}}(t)\ $for $X_{1}$
and $X_{2}$ can be written as: 
\begin{eqnarray}
m_{w_{i}}(t) &=&\frac{1}{F_{X_{i}}(t)}\int_{0}^{t}F_{X_{i}}(x_{i})\ dx_{i}\ 
\notag \\
&=&\frac{1}{F_{X_{i}}(t)}\dsum\limits_{j,k=0}^{\infty }\frac{(-1)^{j}g_{\ast
j}^{(k)}}{(k+1)!}\binom{\theta _{i}+\theta _{3}}{j}t^{k+1};\ \ \ i=1,2,
\label{108}
\end{eqnarray}%
where 
\begin{equation*}
g_{\ast j}^{(k)}=\frac{d^{k}}{dx_{i}^{k}}(e^{-j\eta ^{\alpha
}(x_{i})})|_{x_{i}=0}.
\end{equation*}

\subsection{The hazard rate, the availability and the mean residual lifetime
functions for a parallel system}

Cox (1972) defined the joint hazard rate function as a vector, which is
useful to calculate the total life span of a two component parallel system $%
(2-out-of-2:F)$ as follows%
\begin{equation}
h(x^{\ast })=\left( h_{X}(x),h_{12}(x_{1}|x_{2}),h_{21}(x_{2}|x_{1})\right) ,
\label{109}
\end{equation}%
where the first element $h(x)\ $in the vector $h(x^{\ast })$, gives the
hazard function of the system using the information that both the component
has survived beyond $x$, where$\ X=\min (X_{1},X_{2})$. The second element $%
h_{12}(x_{1}|x_{2})$,$\ $gives the hazard function span of the first
component given that it has survived to an age $x_{1}$, and the other has
failed at $x_{2}$. Similar argument holds for the third element $%
h_{21}(x_{2}|x_{1})$.

If ($X_{1},X_{2}$) is a BEGLE random vector, then the joint hazard rate
function $h(x^{\ast })$ is

\begin{equation}
h_{X}(x)=\frac{\alpha (\theta _{1}+\theta _{2}+\theta _{3})(a+bx)\eta
^{\alpha -1}(x)e^{-\eta ^{\alpha }(x)}\left( 1-e^{-\eta ^{\alpha
}(x)}\right) ^{-1}}{\left( 1-e^{-\eta ^{\alpha }(x)}\right) ^{\theta
_{1}+\theta _{2}+\theta _{3}-1}-1},
\end{equation}

\begin{equation}
h_{12}(x_{1}|x_{2})=\frac{f_{X_{1}}(x_{1})\left( 1-e^{-\eta ^{\alpha
}(x_{1})}\right) ^{\theta _{3}}}{1-\left( 1-e^{-\eta ^{\alpha
}(x_{1})}\right) ^{\theta _{1}}}
\end{equation}%
and

\begin{equation}
h_{21}(x_{2}|x_{1})=\frac{f_{X_{2}}(x_{2})\left( 1-e^{-\eta ^{\alpha
}(x_{2})}\right) ^{\theta _{3}}}{1-\left( 1-e^{-\eta ^{\alpha
}(x_{2})}\right) ^{\theta _{2}}}.
\end{equation}%
Also, the joint availability function can be defined as a vector, which is
useful to calculate the expected lifetime of a parallel system with two
component as follows 
\begin{equation}
V(x^{\ast })=\left( v_{X}(x),v_{12}(x_{1}|x_{2}),v_{21}(x_{2}|x_{1})\right) ,
\label{110}
\end{equation}%
where the first element $v(x)\ $in the vector $V(x^{\ast })$,$\ $gives the
expected lifetime of the system using the information that both the
component has survived beyond $x$,$\ $where$\ X=\min (X_{1},X_{2})$. The
second element $v_{12}(x_{1}|x_{2})$,$\ $gives the expected lifetime span of
the first component given that it has survived to an age $x_{1}$, and the
other has failed at $x_{2}$. Similar argument holds for the third element $%
v_{21}(x_{2}|x_{1})$.

If ($X_{1},X_{2}$) is a BEGLE random vector, then the joint availability
function$V(x^{\ast })$ is

\begin{equation}
V(x)=\frac{1}{A}\dint\limits_{x}^{\infty }yf_{X}(y)\ dy;\
A=\dint\limits_{x}^{\infty }f_{X}(y)\ dy\ ,
\end{equation}%
using Maclaurin expansion, binomial expansion and upper incomplete gamma
function, we get 
\begin{equation}
V(x)=\frac{1}{1-\left( \Psi (x)\right) ^{\theta _{1}+\theta _{2}+\theta _{3}}%
}\sum_{i=0}^{\infty }\sum_{j=0}^{\alpha -1}\sum_{k=0}^{\infty }\vartheta
_{i,j}^{(k)}(a\Gamma (\alpha ^{\ast },x)+b(\Gamma (\alpha ^{\ast }+1,x)).
\end{equation}%
Similarly,

\begin{eqnarray}
V_{12}(x_{1}|x_{2}) &=&\frac{1}{B}\dint\limits_{x_{1}}^{\infty }yf(y,x_{2})\
dy;\ B=\dint\limits_{x_{1}}^{\infty }f(y,x_{2})\ dy,\ x_{1}>x_{2},  \notag \\
&=&\frac{1}{1-\left( \Psi (x_{1})\right) ^{\theta _{1}}}\sum_{j=0}^{\alpha
-1}\sum_{i,k=0}^{\infty }\xi _{i,j}^{(k)}(a\Gamma (\alpha ^{\ast
},x_{1})+b(\Gamma (\alpha ^{\ast }+1,x_{1})),
\end{eqnarray}%
and

\begin{eqnarray}
V_{21}(x_{2}|x_{1}) &=&\frac{1}{C}\dint\limits_{x_{2}}^{\infty }yf(x_{1},y)\
dy;\ C=\dint\limits_{x_{2}}^{\infty }f(x_{1},y)\ dy,\ x_{1}<x_{2},  \notag \\
&=&\frac{1}{1-\left( \Psi (x_{2})\right) ^{\theta _{2}}}\sum_{j=0}^{\alpha
-1}\sum_{i,k=0}^{\infty }\Omega _{i,j}^{(k)}(a\Gamma (\alpha ^{\ast
},x_{2})+b(\Gamma (\alpha ^{\ast }+1,x_{2})),
\end{eqnarray}%
where

\begin{eqnarray*}
\vartheta _{i,j}^{(k)} &=&\alpha \left( \theta _{1}+\theta _{2}+\theta
_{3}\right) \frac{(-1)^{i}a^{\alpha -j-1}Q_{\ast i}^{(k)}}{k!}\left( \frac{b%
}{2}\right) ^{j}\binom{\alpha -1}{j}\binom{\theta _{1}+\theta _{2}+\theta
_{3}-1}{i}, \\
\xi _{i,j}^{(k)} &=&\alpha \theta _{1}\frac{(-1)^{i}a^{\alpha -j-1}Q_{\ast
i}^{(k)}}{k!}\left( \frac{b}{2}\right) ^{j}\binom{\alpha -1}{j}\binom{\theta
_{1}-1}{i}, \\
\Omega _{i,j}^{(k)} &=&\alpha \theta _{2}\frac{(-1)^{i}a^{\alpha
-j-1}Q_{\ast i}^{(k)}}{k!}\left( \frac{b}{2}\right) ^{j}\binom{\alpha -1}{j}%
\binom{\theta _{2}-1}{i}, \\
\alpha ^{\ast } &=&\alpha +j+k+1,
\end{eqnarray*}%
and 
\begin{equation*}
Q_{\ast i}^{(k)}=\frac{d^{k}}{dy^{k}}(\exp (-(1+i)\eta ^{\alpha
}(y)+y))|_{y=0}.
\end{equation*}%
On the other hand, Asha and Jagathnath (2008) defined the joint mean
residual lifetime $m(x^{\ast })$, which is useful to compute the mean
residual lifetime (MRL) to two component in a parallel system, as follows 
\begin{equation}
m(x^{\ast })=\left( m_{X}(x),m_{12}(x_{1}|x_{2}),m_{21}(x_{2}|x_{1})\right) ,
\end{equation}%
where the first element $m(x)\ $in the vector $m(x^{\ast })$,$\ $gives the
MRL of the system using the information that both the component has survived
beyond $x\ $,where$\ X=\min (X_{1},X_{2})$. The second element $%
m_{12}(x_{1}|x_{2})$,$\ $gives the MRL span of the first component given
that it has survived to an age $x_{1}$, and the other has failed at $x_{2}$.
Similar argument holds for the third element $m_{21}(x_{2}|x_{1})$. The
joint MRL function related to the joint vitality function by the
relationships%
\begin{eqnarray}
m_{X}(x) &=&v_{X}(x)-x\ \ \ ,x>0. \\
m_{12}(x_{1}|x_{2}) &=&v_{12}(x_{1}|x_{2})-x_{1}\ \ ,\ \ x_{1}>x_{2}. \\
m_{21}(x_{2}|x_{1}) &=&v_{21}(x_{2}|x_{1})-x_{2}\ \ ,\ \ x_{1}<x_{2}.
\end{eqnarray}%
So, If ($X_{1},X_{2}$) is a BEGLE random vector, then it is easy to get the
vector $m(x^{\ast })$.

\section{Maximum Likelihood Estimation (MLE)}

In this section, we want to estimate the unknown parameters of the BEGLED.
We will use the maximum likelihood method. Suppose that $(x_{11},x_{21})$, $%
(x_{12},x_{22})$,..., $(x_{1n},x_{2n})$ is a sample of size n, from the
BEGLED. We use the following notation $I_{1}=\{x_{1i}<x_{2i}\}$,\ $%
I_{2}=\{x_{1i}>x_{2i}\}$,\ $I_{3}=\{x_{1i}=x_{2i}=x_{i}\}$,$\ I=I_{1}\cup
I_{2}\cup I_{3}$, $\left\vert I_{1}\right\vert =n_{1},$ $\left\vert
I_{2}\right\vert =n_{2},$ $\ \left\vert I_{3}\right\vert =n_{3},$ and $%
\left\vert I\right\vert =n_{1}+n_{2}+n_{3}=n.$ Based on the observations,
the likelihood function $l(\Phi )\ $of this sample is

\begin{equation}
l(\Phi
)=\dprod\limits_{i=1}^{n_{1}}f_{1}(x_{1i},x_{2i})\dprod%
\limits_{i=1}^{n_{2}}f_{2}(x_{1i},x_{2i})\dprod%
\limits_{i=1}^{n_{3}}f_{3}(x_{i},x_{i}).  \label{55}
\end{equation}%
Substituting from Equation (\ref{3}) into Equation (\ref{55}), the
log-likelihood function $L(\Phi )\ $can be written as%
\begin{eqnarray}
L(\Phi ) &=&n_{1}\ln \left( \alpha ^{2}\theta _{2}\left( \theta _{1}+\theta
_{3}\right) \right) +\dsum\limits_{i=1}^{{\small n}_{1}}\ln
(a+bx_{1i})+(\alpha -1)\dsum\limits_{i=1}^{{\small n}_{1}}\ln (\eta (x_{1i}))%
\text{\ }  \notag \\
&&+\left( \theta _{1}+\theta _{3}-1\right) \dsum\limits_{i=1}^{{\small n}%
_{1}}\ln \left( 1-e^{-\eta ^{\alpha }(x_{1i})}\right) -\dsum\limits_{i=1}^{%
{\small n}_{1}}\eta ^{\alpha }(x_{1i})+\dsum\limits_{i=1}^{{\small n}%
_{1}}\ln (a+bx_{2i})  \notag \\
&&+(\alpha -1)\dsum\limits_{i=1}^{{\small n}_{1}}\ln (\eta (x_{2i}))+\left(
\theta _{2}-1\right) \ln \left( 1-e^{-\eta ^{\alpha }(x_{2i})}\right)
-\dsum\limits_{i=1}^{{\small n}_{1}}\eta ^{\alpha }(x_{2i})  \notag \\
&&+n_{2}\ln (\alpha ^{2}\theta _{1}(\theta _{2}+\theta
_{3}))+\dsum\limits_{i=1}^{{\small n}_{2}}\ln (a+bx_{1i})+(\alpha
-1)\dsum\limits_{i=1}^{{\small n}_{2}}\ln (\eta (x_{1i}))  \notag \\
&&+(\theta _{1}-1)\dsum\limits_{i=1}^{{\small n}_{2}}\ln \left( 1-e^{-\eta
^{\alpha }(x_{1i})}\right) -\dsum\limits_{i=1}^{{\small n}_{2}}\eta ^{\alpha
}(x_{1i})+\dsum\limits_{i=1}^{{\small n}_{2}}\ln (a+bx_{2i})  \notag \\
&&+(\alpha -1)\dsum\limits_{i=1}^{{\small n}_{2}}\ln (\eta (x_{2i}))+\left(
\theta _{2}+\theta _{3}-1\right) \dsum\limits_{i=1}^{{\small n}_{2}}\ln
\left( 1-e^{-\eta ^{\alpha }(x_{2i})}\right)  \notag \\
&&-\dsum\limits_{i=1}^{{\small n}_{2}}\eta ^{\alpha }(x_{2i})+n_{3}\ln 
{\small \alpha }\theta _{3}+\dsum\limits_{i=1}^{{\small n}_{3}}\ln
(a+bx_{i})+(\alpha -1)\dsum\limits_{i=1}^{{\small n}_{3}}\ln (\eta (x_{i})) 
\notag \\
&&+\left( \theta _{1}+\theta _{2}+\theta _{3}-1\right) \dsum\limits_{i=1}^{%
{\small n}_{3}}\ln \left( 1-e^{-\eta ^{\alpha }(x_{i})}\right)
-\dsum\limits_{i=1}^{{\small n}_{3}}\eta ^{\alpha }(x_{i}).  \label{70}
\end{eqnarray}%
The first partial derivatives of Equation (\ref{70}) with respect to $\alpha
,a,b,\theta _{1},\theta _{2}$ and $\theta _{3}$ are

\begin{eqnarray}
\frac{\partial L}{\partial \alpha } &{\small =}&\frac{2n_{1}}{\alpha }%
+\dsum\limits_{i=1}^{{\small n}_{1}}\ln (\eta (x_{1i}))+(\theta _{1}+\theta
_{3}-1)\dsum\limits_{i=1}^{{\small n}_{1}}\frac{\eta ^{\alpha }(x_{1i})\ln
(\eta (x_{1i}))}{e^{\eta ^{\alpha }(x_{1i})}-1}  \notag \\
&&-\dsum\limits_{i=1}^{{\small n}_{1}}\eta ^{\alpha }(x_{1i})\ln (\eta
(x_{1i}))+\dsum\limits_{i=1}^{{\small n}_{1}}\ln (\eta
(x_{2i}))-\dsum\limits_{i=1}^{{\small n}_{1}}\eta ^{\alpha }(x_{2i})\ln
(\eta (x_{2i}))  \notag \\
&&+(\theta _{2}-1)\dsum\limits_{i=1}^{{\small n}_{1}}\frac{\eta ^{\alpha
}(x_{2i})\ln (\eta (x_{2i}))}{e^{\eta ^{\alpha }(x_{2i})}-1}+\frac{2n_{2}}{%
\alpha }+\dsum\limits_{i=1}^{{\small n}_{2}}\ln (\eta (x_{1i}))  \notag \\
&&+(\theta _{1}-1)\dsum\limits_{i=1}^{{\small n}_{2}}\frac{\eta ^{\alpha
}(x_{1i})\ln (\eta (x_{1i}))}{e^{\eta ^{\alpha }(x_{1i})}-1}%
-\dsum\limits_{i=1}^{{\small n}_{2}}\eta ^{\alpha }(x_{1i})\ln (\eta
(x_{1i}))  \notag \\
&&+\dsum\limits_{i=1}^{{\small n}_{2}}\ln (\eta (x_{1i}))+(\theta
_{2}+\theta _{3}-1)\dsum\limits_{i=1}^{{\small n}_{2}}\frac{\eta ^{\alpha
}(x_{2i})\ln (\eta (x_{2i}))}{e^{\eta ^{\alpha }(x_{2i})}-1}  \notag \\
&&-\dsum\limits_{i=1}^{{\small n}_{2}}\eta ^{\alpha }(x_{2i})\ln (\eta
(x_{2i}))+\dsum\limits_{i=1}^{{\small n}_{3}}\ln (\eta
(x_{i}))-\dsum\limits_{i=1}^{{\small n}_{3}}\eta ^{\alpha }(x_{i})\ln (\eta
(x_{i}))  \notag \\
&&+\frac{n_{3}}{\alpha }+(\theta _{1}+\theta _{2}+\theta
_{3}-1)\dsum\limits_{i=1}^{{\small n}_{3}}\frac{\eta ^{\alpha }(x_{i})\ln
(\eta (x_{i}))}{e^{\eta ^{\alpha }(x_{i})}-1},
\end{eqnarray}%
\begin{eqnarray}
\frac{\partial L}{\partial a} &{\tiny =}&\dsum\limits_{i=1}^{{\small n}_{1}}%
\frac{1}{a+bx_{1i}}+(\alpha -1)\dsum\limits_{i=1}^{{\small n}_{1}}\frac{%
x_{1i}}{\eta (x_{1i})}+(\theta _{1}+\theta _{3}-1)\dsum\limits_{i=1}^{%
{\small n}_{1}}\frac{\alpha x_{1i}\eta ^{\alpha -1}(x_{1i})}{e^{\eta
^{\alpha }(x_{1i})}-1}  \notag \\
&&-\alpha \dsum\limits_{i=1}^{{\small n}_{1}}x_{1i}\eta ^{\alpha
-1}(x_{1i})+\dsum\limits_{i=1}^{{\small n}_{1}}\frac{1}{a+bx_{2i}}+(\alpha
-1)\dsum\limits_{i=1}^{{\small n}_{1}}\frac{x_{2i}}{\eta (x_{2i})}  \notag \\
&&+(\theta _{2}-1)\dsum\limits_{i=1}^{{\small n}_{1}}\frac{\alpha x_{2i}\eta
^{\alpha -1}(x_{2i})}{e^{\eta ^{\alpha }(x_{2i})}-1}-\alpha
\dsum\limits_{i=1}^{{\small n}_{1}}x_{2i}\eta ^{\alpha
-1}(x_{2i})+\dsum\limits_{i=1}^{{\small n}_{2}}\frac{1}{a+bx_{1i}}  \notag \\
&&+(\alpha -1)\dsum\limits_{i=1}^{{\small n}_{2}}\frac{x_{1i}}{\eta (x_{1i})}%
-\alpha \dsum\limits_{i=1}^{{\small n}_{2}}x_{1i}\eta ^{\alpha
-1}(x_{1i})+(\theta _{1}-1)\dsum\limits_{i=1}^{{\small n}_{2}}\frac{\alpha
x_{1i}\eta ^{\alpha -1}(x_{1i})}{e^{\eta ^{\alpha }(x_{1i})}-1}  \notag \\
&&+\dsum\limits_{i=1}^{{\small n}_{2}}\frac{1}{a+bx_{2i}}+(\alpha
-1)\dsum\limits_{i=1}^{{\small n}_{2}}\frac{x_{2i}}{\eta (x_{2i})}+(\theta
_{2}+\theta _{3}-1)\dsum\limits_{i=1}^{{\small n}_{2}}\frac{\alpha
x_{2i}\eta ^{\alpha -1}(x_{2i})}{e^{\eta ^{\alpha }(x_{2i})}-1}  \notag \\
&&-\alpha \dsum\limits_{i=1}^{{\small n}_{2}}x_{2i}\eta ^{\alpha
-1}(x_{2i})+\dsum\limits_{i=1}^{{\small n}_{3}}\frac{1}{a+bx_{i}}+(\theta
_{1}+\theta _{2}+\theta _{3}-1)\dsum\limits_{i=1}^{{\small n}_{3}}\frac{%
\alpha x_{i}\eta ^{\alpha -1}(x_{i})}{e^{\eta ^{\alpha }(x_{i})}-1}  \notag
\\
&&+(\alpha -1)\dsum\limits_{i=1}^{{\small n}_{3}}\frac{x_{i}}{\eta (x_{i})}%
-\alpha \dsum\limits_{i=1}^{{\small n}_{3}}x_{i}\eta ^{\alpha -1}(x_{i}),
\end{eqnarray}%
\begin{eqnarray}
\frac{\partial L}{\partial b} &{\tiny =}&\dsum\limits_{i=1}^{{\small n}_{1}}%
\frac{x_{1i}}{a+bx_{1i}}+\frac{\alpha -1}{2}\dsum\limits_{i=1}^{{\small n}%
_{1}}\frac{(x_{1i})^{2}}{\eta (x_{1i})}+\frac{\theta _{1}+\theta _{3}-1}{2}%
\dsum\limits_{i=1}^{{\small n}_{1}}\frac{\alpha (x_{1i})^{2}\eta ^{\alpha
-1}(x_{1i})}{e^{\eta ^{\alpha }(x_{1i})}-1}  \notag \\
&&+\dsum\limits_{i=1}^{{\small n}_{1}}\frac{x_{2i}}{a+bx_{2i}}-\frac{\alpha 
}{2}\dsum\limits_{i=1}^{{\small n}_{1}}(x_{1i})^{2}\eta ^{\alpha -1}(x_{1i})+%
\frac{\theta _{2}-1}{2}\dsum\limits_{i=1}^{{\small n}_{1}}\frac{\alpha
(x_{2i})^{2}\eta ^{\alpha -1}(x_{2i})}{e^{\eta ^{\alpha }(x_{2i})}-1}  \notag
\\
&&+\frac{\alpha -1}{2}\dsum\limits_{i=1}^{{\small n}_{1}}\frac{(x_{2i})^{2}}{%
\eta (x_{2i})}-\frac{\alpha }{2}\dsum\limits_{i=1}^{{\small n}%
_{1}}(x_{2i})^{2}\eta ^{\alpha -1}(x_{2i})+\dsum\limits_{i=1}^{{\small n}%
_{2}}\frac{x_{1i}}{a+bx_{1i}}  \notag \\
&&+\frac{\alpha -1}{2}\dsum\limits_{i=1}^{{\small n}_{2}}\frac{(x_{1i})^{2}}{%
\eta (x_{1i})}+\frac{\theta _{1}-1}{2}\dsum\limits_{i=1}^{{\small n}_{2}}%
\frac{\alpha (x_{1i})^{2}\eta ^{\alpha -1}(x_{1i})}{e^{\eta ^{\alpha
}(x_{1i})}-1}  \notag \\
&&-\frac{\alpha }{2}\dsum\limits_{i=1}^{{\small n}_{2}}(x_{1i})^{2}\eta
^{\alpha -1}(x_{1i})+\dsum\limits_{i=1}^{{\small n}_{2}}\frac{x_{2i}}{%
a+bx_{2i}}+\frac{\alpha -1}{2}\dsum\limits_{i=1}^{{\small n}_{2}}\frac{%
(x_{2i})^{2}}{\eta (x_{2i})}  \notag \\
&&+\frac{\theta _{2}+\theta _{3}-1}{2}\dsum\limits_{i=1}^{{\small n}_{2}}%
\frac{\alpha (x_{2i})^{2}\eta ^{\alpha -1}(x_{2i})}{e^{\eta ^{\alpha
}(x_{2i})}-1}-\frac{\alpha }{2}\dsum\limits_{i=1}^{{\small n}%
_{2}}(x_{2i})^{2}\eta ^{\alpha -1}(x_{2i})  \notag \\
&&+\dsum\limits_{i=1}^{{\small n}_{3}}\frac{x_{i}}{a+bx_{i}}+\frac{\alpha -1%
}{2}\dsum\limits_{i=1}^{{\small n}_{3}}\frac{(x_{i})^{2}}{\eta (x_{i})}+%
\frac{\theta _{1}+\theta _{2}+\theta _{3}-1}{2}\dsum\limits_{i=1}^{{\small n}%
_{3}}\frac{\alpha (x_{i})^{2}\eta ^{\alpha -1}(x_{i})}{e^{\eta ^{\alpha
}(x_{i})}-1}  \notag \\
&&-\frac{\alpha }{2}\dsum\limits_{i=1}^{{\small n}_{3}}(x_{i})^{2}\eta
^{\alpha -1}(x_{i}),
\end{eqnarray}%
\begin{eqnarray}
\frac{\partial L}{\partial \theta _{1}} &{\tiny =}&\frac{n_{1}}{\theta
_{1}+\theta _{3}}+\dsum\limits_{i=1}^{{\small n}_{1}}\ln \left( 1-e^{-\eta
^{\alpha }(x_{1i})}\right) +\frac{n_{2}}{\theta _{1}}+\dsum\limits_{i=1}^{%
{\small n}_{2}}\ln \left( 1-e^{-\eta ^{\alpha }(x_{1i})}\right)  \notag \\
&&+\dsum\limits_{i=1}^{{\small n}_{3}}\ln \left( 1-e^{-\eta ^{\alpha
}(x_{i})}\right) ,
\end{eqnarray}%
\begin{eqnarray}
\frac{\partial L}{\partial \theta _{2}} &{\tiny =}&\frac{n_{1}}{\theta _{2}}%
+\dsum\limits_{i=1}^{{\small n}_{1}}\ln \left( 1-e^{-\eta ^{\alpha
}(x_{2i})}\right) +\frac{n_{2}}{\theta _{2}+\theta _{3}}+\dsum\limits_{i=1}^{%
{\small n}_{2}}\ln \left( 1-e^{-\eta ^{\alpha }(x_{2i})}\right)  \notag \\
&&+\dsum\limits_{i=1}^{{\small n}_{3}}\ln \left( 1-e^{-\eta ^{\alpha
}(x_{i})}\right) ,
\end{eqnarray}%
and%
\begin{eqnarray}
\frac{\partial L}{\partial \theta _{3}} &{\tiny =}&\frac{n_{1}}{\theta
_{1}+\theta _{3}}+\dsum\limits_{i=1}^{{\small n}_{1}}\ln \left( 1-e^{-\eta
^{\alpha }(x_{1i})}\right) +\dsum\limits_{i=1}^{{\small n}_{2}}\ln \left(
1-e^{-\eta ^{\alpha }(x_{2i})}\right) +\frac{n_{3}}{\theta _{3}}  \notag \\
&&+\frac{n_{2}}{\theta _{2}+\theta _{3}}+\dsum\limits_{i=1}^{{\small n}%
_{3}}\ln \left( 1-e^{-\eta ^{\alpha }(x_{i})}\right) .
\end{eqnarray}%
By Equating the Equations (42-47) by zeros, we get the non-linear normal
Equations. So, the solution has to be obtained numerically.

\section{Data Analysis}

In this section, we have analyzed one bivariate real data set to explicate
that the BEGLED can be a good lifetime model, comparing with Marshall-Olkin
bivariate exponential distribution (MOBED), bivariate generalized
exponential distribution (BVGED) and bivariate generalized linear failure
rate distribution (BGLFRD). To make this comparison, we will use the
log-likelihood values (L), Akaike information criterion (AIC), correct
Akaike information criterion (CAIC), Hannan-Quinn information criterion
(HQIC) and the likelihood ratio test ($\Lambda $).

The data set in Table 1 has been obtained from Meintanis (2007). This data
represents football (soccer) data of the UEFA Champion's League data for the
year 2004 : 2005 and 2005 : 2006. This data describes the games which
satisfy\ the following two conditions:

\begin{enumerate}
\item At least one kick goal scored by any team have been considered.

\item The home team must be scored at least one goal.
\end{enumerate}

Note that, the kick goal is the goal which scored directly from foul kick,
penalty kick or any other direct free kick. Here\ the variables $X_{1}$ and $%
X_{2}$ are as follows:

$X_{1}$ : represents the time in minutes of the first kick goal scored by
any team.

$X_{2}$ : represents the first goal of any type scored by the home team.

\begin{equation*}
\begin{tabular}{|ccc|ccc|}
\multicolumn{6}{l|}{$\text{\textbf{Table 1.} The UEFA Champion's League data
for the year 2004 : 2005 and 2005 : 2006.}$} \\ \hline\hline
$2005:2006$ \ \ \ \ \ \ \ \ \ \ \  & $\ \ \ X_{1}$ \ \ \ \  & $\ \ X_{2}$ \
\  & $2004:2005$ \ \ \ \ \ \ \ \  & $\ \ X_{1}$ \ \  & $\ \ X_{2}$ \ \  \\ 
\hline\hline
\multicolumn{1}{|l}{\small Lyon : Real Madrid} & ${\tiny 26}$ & ${\tiny 20}$
& \multicolumn{1}{|l}{\small Internazionale : Bremen} & ${\tiny 34}$ & $%
{\tiny 34}$ \\ 
\multicolumn{1}{|l}{\small Milan : Fenerbahce} & ${\tiny 63}$ & ${\tiny 18}$
& \multicolumn{1}{|l}{\small Real Madrid : Roma} & ${\tiny 53}$ & ${\tiny 39}
$ \\ 
\multicolumn{1}{|l}{\small Chelsea : Anderlecht} & ${\tiny 19}$ & ${\tiny 19}
$ & \multicolumn{1}{|l}{\small Man. United : Fenerbahce} & ${\tiny 54}$ & $%
{\tiny 7}$ \\ 
\multicolumn{1}{|l}{\small Club Brugge : Juventus} & ${\tiny 66}$ & ${\tiny %
85}$ & \multicolumn{1}{|l}{\small Bayern : Ajax} & ${\tiny 51}$ & ${\tiny 28}
$ \\ 
\multicolumn{1}{|l}{\small Fenerbahce : PSV} & ${\tiny 40}$ & ${\tiny 40}$ & 
\multicolumn{1}{|l}{\small Moscow : PSG} & ${\tiny 76}$ & ${\tiny 64}$ \\ 
\multicolumn{1}{|l}{\small Internazionale : Rangers} & ${\tiny 49}$ & $%
{\tiny 49}$ & \multicolumn{1}{|l}{\small Barcelona : Shakhtar} & ${\tiny 64}$
& ${\tiny 15}$ \\ 
\multicolumn{1}{|l}{\small Panathinaikos : Bremen} & ${\tiny 8}$ & ${\tiny 8}
$ & \multicolumn{1}{|l}{\small Leverkusen : Roma} & ${\tiny 26}$ & ${\tiny 48%
}$ \\ 
\multicolumn{1}{|l}{\small Ajax : Arsenal} & ${\tiny 69}$ & ${\tiny 71}$ & 
\multicolumn{1}{|l}{\small Arsenal :\ Panathinaikos} & ${\tiny 16}$ & $%
{\tiny 16}$ \\ 
\multicolumn{1}{|l}{\small Man. United : Benfica} & ${\tiny 39}$ & ${\tiny 39%
}$ & \multicolumn{1}{|l}{\small Dynamo Kyiv : Real Madrid} & ${\tiny 44}$ & $%
{\tiny 13}$ \\ 
\multicolumn{1}{|l}{\small Real Madrid : Rosenborg} & ${\tiny 82}$ & ${\tiny %
48}$ & \multicolumn{1}{|l}{\small Man. United : Sparta} & ${\tiny 25}$ & $%
{\tiny 14}$ \\ 
\multicolumn{1}{|l}{\small Villarreal : Benfica} & ${\tiny 72}$ & ${\tiny 72}
$ & \multicolumn{1}{|l}{\small Bayern : M. TelAviv} & ${\tiny 55}$ & ${\tiny %
11}$ \\ 
\multicolumn{1}{|l}{\small Juventus : Bayern} & ${\tiny 66}$ & ${\tiny 62}$
& \multicolumn{1}{|l}{\small Bremen : Internazionale} & ${\tiny 49}$ & $%
{\tiny 49}$ \\ 
\multicolumn{1}{|l}{\small Club Brugge : Rapid} & ${\tiny 25}$ & ${\tiny 9}$
& \multicolumn{1}{|l}{\small Anderlecht : Valencia} & ${\tiny 24}$ & ${\tiny %
24}$ \\ 
\multicolumn{1}{|l}{\small Olympiacos : Lyon} & ${\tiny 41}$ & ${\tiny 3}$ & 
\multicolumn{1}{|l}{\small Panathinaikos : PSV} & ${\tiny 44}$ & ${\tiny 30}$
\\ 
\multicolumn{1}{|l}{\small Internazionale : Porto} & ${\tiny 16}$ & ${\tiny %
75}$ & \multicolumn{1}{|l}{\small Arsenal : Rosenborg} & ${\tiny 42}$ & $%
{\tiny 3}$ \\ 
\multicolumn{1}{|l}{\small Schalke : PSV} & ${\tiny 18}$ & ${\tiny 18}$ & 
\multicolumn{1}{|l}{\small Liverpool : Olympiacos} & ${\tiny 27}$ & ${\tiny %
47}$ \\ 
\multicolumn{1}{|l}{\small Barcelona : Bremen} & ${\tiny 22}$ & ${\tiny 14}$
& \multicolumn{1}{|l}{\small M. Tel-Aviv : Juventus} & ${\tiny 28}$ & $%
{\tiny 28}$ \\ 
\multicolumn{1}{|l}{\small Milan : Schalke} & ${\tiny 42}$ & ${\tiny 42}$ & 
\multicolumn{1}{|l}{\small Bremen : Panathinaikos} & ${\tiny 2}$ & ${\tiny 2}
$ \\ 
\multicolumn{1}{|l}{\small Rapid : Juventus} & ${\tiny 36}$ & ${\tiny 52}$ & 
\multicolumn{1}{|l}{} &  &  \\ \hline\hline
\end{tabular}%
\end{equation*}%
\newline
To analyze this data by the BEGLED, we fit at first the marginals $X_{1}$
and $X_{2}$ of the BEGLED separately one by one on this data. The following
Tables obtain the MLEs, L, Anderson-Darling (A$^{\ast }$) and Cram\'{e}r-Von
Mises (W$^{\ast }$)\textbf{\ }values\textbf{\ }for the marginals $X_{1}$ and 
$X_{2}$ respectively for each model. 
\begin{equation*}
\begin{tabular}{cccccccc}
\multicolumn{8}{l}{$\text{\textbf{Table 2.} The MLE(s), L, }A^{\ast }\text{, 
}W^{\ast }\text{ values for }X_{1}.$} \\ \hline\hline
\ \ \ \ {\small Model \ \ } & $\overset{\wedge }{{\tiny a}}$ & $\overset{%
\wedge }{{\tiny b}}$ & {\small \ }$\overset{\wedge }{{\tiny \theta }}$ & $%
\overset{\wedge }{{\tiny \alpha }}$ & {\small -L} & A$^{\ast }$ & W$^{\ast }$
\\ \hline\hline
{\small E} & ${\tiny 0.0245}$ & - & - & - & ${\tiny 174.30}$ & ${\tiny 0.5202%
}$ & ${\tiny 0.0686}$ \\ 
{\small GE} & ${\tiny 0.0449}$ & - & ${\tiny 3.119}$ & - & ${\tiny 165.82}$
& ${\tiny 0.6171}$ & ${\tiny 0.0826}$ \\ 
{\small GLFR} & ${\tiny 0.0052}$ & ${\tiny 0.0009}$ & ${\tiny 1.302}$ & - & $%
{\tiny 162.68}$ & ${\tiny 0.2637}$ & ${\tiny 0.0399}$ \\ 
{\small EGLE} & ${\tiny 0.0022}$ & ${\tiny 0.0006}$ & ${\tiny 0.492}$ & $%
{\tiny 1.897}$ & ${\tiny 161.89}$ & ${\tiny 0.2530}$ & ${\tiny 0.0396}$ \\ 
\hline\hline
\end{tabular}%
\end{equation*}%
\begin{equation*}
\begin{tabular}{cccccccc}
\multicolumn{8}{l}{$\text{\textbf{Table 3.} The MLE(s), L, }A^{\ast }\text{, 
}W^{\ast }\text{ values for }X_{2}.$} \\ \hline\hline
\ \ {\small Model \ \ \ \ } & $\overset{\wedge }{{\tiny a}}$ & $\overset{%
\wedge }{{\tiny b}}$ & $\overset{\wedge }{{\tiny \theta }}$ & $\overset{%
\wedge }{{\tiny \alpha }}$ & -{\small L} & \ A$^{\ast }$ & W$^{\ast }$ \\ 
\hline\hline
{\small E} & ${\tiny 0.0304}$ & {\tiny --} & {\tiny --} & {\tiny --} & $%
{\tiny 166.219}$ & ${\tiny 0.3651}$ & ${\tiny 0.0549}$ \\ 
{\small GE} & ${\tiny 0.0413}$ & {\tiny --} & ${\tiny 1.678}$ & {\tiny --} & 
${\tiny 163.937}$ & ${\tiny 0.3859}$ & ${\tiny 0.0576}$ \\ 
{\small GLFR} & ${\tiny 0.0192}$ & ${\tiny 6\times 10}^{{\tiny -4}}$ & $%
{\tiny 1.14}$ & {\tiny --} & ${\tiny 162.938}$ & ${\tiny 0.2713}$ & ${\tiny %
0.04478}$ \\ 
{\small EGLE} & ${\tiny 0.0172}$ & ${\tiny 2\times 10}^{{\tiny -4}}$ & $%
{\tiny 0.622}$ & ${\tiny 1.705}$ & ${\tiny 162.672}$ & ${\tiny 0.2640}$ & $%
{\tiny 0.0436}$ \\ \hline\hline
\end{tabular}%
\end{equation*}%
We can conclude that, the EGLE distribution fits the data better than E, GE
and GLFR{\small \ }distributions{\small \ }for the marginals, because it has
the smallest value among -L, A$^{\ast }\ $and W$^{\ast }$.

Since, the E, GE and GLFR{\small \ }distributions are special cases from the
EGLE distribution, we perform the following three testing of hypotheses%
\textbf{\ }for\textbf{\ }$X_{1}$ and $X_{2}$ separately:

\begin{enumerate}
\item[ \ \ \ Test 1:] $H_{01}:$ $\alpha =1,b=0,\theta =1$ (ED) against $%
H_{11}:$ $\alpha \neq 1,b>0,\theta \neq 1$ (EGLED).

\item[ \ \ \ Test 2:] $H_{02}:$ $\alpha =1,b=0$ (GED) against $H_{12}:$ $%
\alpha \neq 1,b>0$ (EGLED).

\item[ \ \ \ Test 3:] $H_{03}:$ $\alpha =1$ (GLFRD) against $H_{13}:$ $%
\alpha \neq 1$ (EGLED).
\end{enumerate}

\bigskip The likelihood ratio test statistics ($\Lambda $), the degree of
freedom (d.f) and the corresponding p-values for the three tests of
hypotheses in case of\textbf{\ }$X_{1}$ and $X_{2}$ are presented in Tables
4 and 5 respectively.

\begin{equation*}
\begin{tabular}{ccccc}
\multicolumn{5}{l}{$\text{\textbf{Table 4.} The likelihood ratio test
statistics, d.f and p-values for }X_{1}.$} \\ \hline\hline
{\small \ \ \ Model \ \ } & {\small \ }$\ \ \ H_{\circ }${\small \ \ \ \ \ \ 
} & {\small \ }$\ \ \ \ \ \ \Lambda ${\small \ \ \ \ \ \ \ } & {\small \ \
d.f. \ \ } & {\small \ \ \ \ \ p-values \ \ \ \ } \\ \hline\hline
{\small E} & ${\tiny \alpha =1,b=0,\theta =1}$ & ${\tiny 24.824}$ & ${\tiny 3%
}$ & ${\tiny 0.00001681}$ \\ 
{\small GE} & ${\tiny \alpha =1,b=0}$ & ${\tiny 7.846}$ & ${\tiny 2}$ & $%
{\tiny 0.01978166}$ \\ 
{\small GLFR} & ${\tiny \alpha =1}$ & ${\tiny 1.576}$ & ${\tiny 1}$ & $%
{\tiny 0.20933780}$ \\ \hline\hline
\end{tabular}%
\end{equation*}

\begin{equation*}
\begin{tabular}{ccccc}
\multicolumn{5}{l}{$\text{\textbf{Table 5.} The likelihood ratio test
statistics, d.f and p-values for }X_{2}.$} \\ \hline\hline
{\small \ \ \ Model \ \ } & {\small \ }$\ \ \ \ \ H_{\circ }${\small \ \ \ \
\ \ \ \ } & {\small \ }$\ \ \ \ \ \ \ \Lambda ${\small \ \ \ \ \ \ \ \ } & 
{\small \ \ \ \ d.f. \ \ \ \ } & {\small \ \ \ \ p-values \ \ \ } \\ 
\hline\hline
{\small E} & ${\tiny \alpha =1,b=0,\theta =1}$ & ${\tiny 7.094}$ & ${\tiny 3}
$ & ${\tiny 0.06896126}$ \\ 
{\small GE} & ${\tiny \alpha =1,b=0}$ & ${\tiny 2.53}$ & ${\tiny 2}$ & $%
{\tiny 0.2822393}$ \\ 
{\small GLFR} & ${\tiny \alpha =1}$ & ${\tiny 0.532}$ & ${\tiny 1}$ & $%
{\tiny 0.46576723}$ \\ \hline\hline
\end{tabular}%
\end{equation*}

When the level of significance $\delta \ $equals 0.05, it is clear that:

\begin{enumerate}
\item[(a)] The EGLED provides a significantly better fit in case of $X_{1}$
and $X_{2}$ compared to the ED.

\item[(b)] The EGLED provides a significantly better fit in case of $X_{1}$
compared to the GED.

\item[(c)] The EGLED provides a better fit for $X_{2}$ compared to the GED.

\item[(d)] The EGLED provides a better fit in case of $X_{1}$ and $X_{2}$
compared to the GLFRD.
\end{enumerate}

On the other hand, after studying the marginals $X_{1}$ and $X_{2},$ we fit
the BEGLED on the UEFA Champion's League data. The following tables obtain
the MLEs, L, AIC, CAIC and HQIC values.

\begin{equation*}
\begin{tabular}{cccccccc}
\multicolumn{8}{l}{$\text{\textbf{Table 6.} The MLEs and L values}$} \\ 
\hline\hline
{\small Model} & $\overset{\wedge }{{\tiny \alpha }}$ & $\overset{\wedge }{%
{\tiny a}}$ & $\overset{\wedge }{{\tiny b}}$ & {\small \ }$\overset{\wedge }{%
{\tiny \theta }_{1}}$ & $\overset{\wedge }{{\tiny \theta }_{2}}$ & $\overset{%
\wedge }{{\tiny \theta }_{3}}$ & {\small -L} \\ \hline\hline
{\small MOBE} & {\tiny --\ } & {\tiny --} & {\tiny --} & ${\tiny 0.012}$ & $%
{\tiny 0.014}$ & ${\tiny 0.022}$ & ${\tiny 339.0}$ \\ 
{\small BVGE} & {\tiny --} & ${\tiny 0.039}$ & {\tiny --} & ${\tiny 1.351}$
& ${\tiny 0.465}$ & ${\tiny 1.153}$ & ${\tiny 296.9}$ \\ 
{\small BGLFR} & {\tiny \ --} & ${\tiny 0.0002}$ & ${\tiny 0.0008}$ & $%
{\tiny 0.492}$ & ${\tiny 0.411}$ & ${\tiny 0.411}$ & ${\tiny 293.4}$ \\ 
{\small BEGLE} & ${\tiny 0.089}$ & ${\tiny 0.0107}$ & ${\tiny 2.711}$ & $%
{\tiny 0.00017}$ & ${\tiny 0.249}$ & ${\tiny 0.220}$ & ${\tiny 291.7}$ \\ 
\hline\hline
\end{tabular}%
\end{equation*}

\begin{equation*}
\begin{tabular}{cccc}
\multicolumn{4}{l}{$\text{\textbf{Table 7.} The AIC, CAIC and HQIC values.}$}
\\ \hline\hline
\ \ \ \ {\small Model \ \ \ \ \ \ \ \ } & \ \ \ \ {\small AIC \ \ } & 
{\small CAIC \ \ \ \ } & {\small HQIC} \\ \hline\hline
{\small MOBE \ } & ${\tiny 684.0}$ & ${\tiny 684.7}$ & ${\tiny 685.8}$ \\ 
{\small BVGE} & ${\tiny 601.9}$ & ${\tiny 603.1}$ & ${\tiny 604.1}$ \\ 
{\small BGLFR} & ${\tiny 596.8}$ & ${\tiny 598.7}$ & ${\tiny 599.6}$ \\ 
{\small BEGLE} & ${\tiny 595.4}$ & ${\tiny 598.2}$ & ${\tiny 598.8}$ \\ 
\hline\hline
\end{tabular}%
\end{equation*}%
It is clear that, the BEGLED provides a better fit than MOBE, BVGE and BGLFR
distributions because it has the smallest value among -L, AIC, CAIC and
HQIC. Since, the BVGE and BGLFR\ distributions are special cases from the
BEGLED, then we perform the following two testing of hypotheses:

\begin{enumerate}
\item[ \ \ \ Test 1:] $H_{02}:$ $\alpha =1,b=0$ (BVGED) against $H_{12}:$ $%
\alpha \neq 1,b>0$ (BEGLED).

\item[ \ \ \ Test 2:] $H_{03}:$ $\alpha =1$ (BVGLFRD) against $H_{13}:$ $%
\alpha \neq 1$ (BEGLED).
\end{enumerate}

The likelihood ratio test statistics, d.f and p-values for the BVGE and the
BGLFR distributions are given in the following Table.

\begin{equation*}
\begin{tabular}{ccccc}
\multicolumn{5}{l}{$\text{\textbf{Table 8.} The likelihood ratio test
statistics, d.f and p-values.}$} \\ \hline\hline
{\small \ \ \ Model \ \ } & {\small \ }$\ \ \ H_{\circ }${\small \ \ \ \ \ \ 
} & {\small \ }$\ \ \ \ \ \ \Lambda ${\small \ \ \ \ \ \ \ } & {\small \ \
d.f. \ \ } & {\small \ \ \ \ \ p-values \ \ \ \ } \\ \hline\hline
{\small BVGE} & $\alpha {\small =1,b=0}$ & ${\small 10.466}$ & ${\small 2}$
& ${\small 0.00533749}$ \\ 
{\small BGLFR} & $\alpha {\small =1}$ & ${\small 3.354}$ & ${\small 1}$ & $%
{\small 0.06704192}$ \\ \hline\hline
\end{tabular}%
\end{equation*}%
We note that the p-value is not large. So, we prefer the BEGLED for
analyzing this data.

\section{Simulation Study}

In this section, the MLE method is used to estimate the parameters $\alpha
,a,b,\theta _{1},\theta _{2}$ and $\theta _{3}$ of the BEGLED. The
population parameters are generated using software "Mathcad prime 3"
package. The sampling distributions are obtained for different sample sizes $%
n=[30,50,100,200]$ from $N=1000$ replications. This study presents an
assessment of the properties of the MLE for the parameters in terms of bias,
variance (Var), mean square error (MSE) and 95\% confidence intervals (C.I),
which be obtained in the following Tables.%
\begin{equation*}
\begin{tabular}{|c|c|c|c|c|c|c|}
\multicolumn{7}{l}{$\text{\textbf{Table 9.} The MLEs, Bias, Var, MSE and C.I
values}$} \\ \hline\hline
n & parameter & Estimate & Bias & Var & MSE & \ C.I \\ \hline\hline
30 & ${\tiny \alpha =1.5}$ & ${\tiny 1.529447}$ & ${\tiny 0.029447}$ & $%
{\tiny 0.0422195}$ & ${\tiny 0.0430866}$ & ${\tiny (1.2060,2.0056)}$ \\ 
\hline
& ${\tiny a=0.5}$ & ${\tiny 0.505473}$ & ${\tiny 0.005473}$ & ${\tiny %
0.0050818}$ & ${\tiny 0.0051117}$ & ${\tiny (0.3717,0.6442)}$ \\ \hline
& ${\tiny b=0.7}$ & ${\tiny 0.711647}$ & ${\tiny 0.011647}$ & ${\tiny %
0.0098601}$ & ${\tiny 0.0099958}$ & ${\tiny (0.5353,0.9118)}$ \\ \hline
& ${\tiny \theta }_{{\tiny 1}}{\tiny =0.8}$ & ${\tiny 0.913251}$ & ${\tiny %
0.113251}$ & ${\tiny 0.2811745}$ & ${\tiny 0.2940002}$ & ${\tiny %
(0.2862,2.2053)}$ \\ \hline
& ${\tiny \theta }_{{\tiny 2}}{\tiny =1.2}$ & ${\tiny 1.300607}$ & ${\tiny %
0.100607}$ & ${\tiny 0.3690383}$ & ${\tiny 0.37916}$ & ${\tiny %
(0.5189,2.7853)}$ \\ \hline
& ${\tiny \theta }_{{\tiny 3}}{\tiny =1.3}$ & ${\tiny 1.35555}$ & ${\tiny %
0.05555}$ & ${\tiny 0.2409224}$ & ${\tiny 0.2440082}$ & ${\tiny %
(0.6489,2.4833)}$ \\ \hline\hline
50 & ${\tiny \alpha =1.5}$ & ${\tiny 1.5125455}$ & ${\tiny 0.0125455}$ & $%
{\tiny 0.0249091}$ & ${\tiny 0.0250665}$ & ${\tiny (1.2673,1.8829)}$ \\ 
\hline
& ${\tiny a=0.5}$ & ${\tiny 0.4996033}$ & ${\tiny -0.000396}$ & ${\tiny %
0.0031267}$ & ${\tiny 0.0031268}$ & ${\tiny (0.3929,0.6071)}$ \\ \hline
& ${\tiny b=0.7}$ & ${\tiny 0.7016369}$ & ${\tiny 0.0016369}$ & ${\tiny %
0.0059615}$ & ${\tiny 0.0059642}$ & ${\tiny (0.5624,0.8559)}$ \\ \hline
& ${\tiny \theta }_{{\tiny 1}}{\tiny =0.8}$ & ${\tiny 0.901861}$ & ${\tiny %
0.101861}$ & ${\tiny 0.1651336}$ & ${\tiny 0.1755093}$ & ${\tiny %
(0.3843,1.9003)}$ \\ \hline
& ${\tiny \theta }_{{\tiny 2}}{\tiny =1.2}$ & ${\tiny 1.2920736}$ & ${\tiny %
0.0920736}$ & ${\tiny 0.2386965}$ & ${\tiny 0.247174}$ & ${\tiny %
(0.5865,2.4482)}$ \\ \hline
& ${\tiny \theta }_{{\tiny 3}}{\tiny =1.3}$ & ${\tiny 1.3579218}$ & ${\tiny %
0.0579218}$ & ${\tiny 0.1491417}$ & ${\tiny 0.1524966}$ & ${\tiny %
(0.7602,2.2375)}$ \\ \hline\hline
100 & ${\tiny \alpha =1.5}$ & ${\tiny 1.5092043}$ & ${\tiny 0.0092043}$ & $%
{\tiny 0.011566}$ & ${\tiny 0.0116508}$ & ${\tiny (1.3322,1.7394)}$ \\ \hline
& ${\tiny a=0.5}$ & ${\tiny 0.5002042}$ & ${\tiny 0.0002042}$ & ${\tiny %
0.0014805}$ & ${\tiny 0.0014805}$ & ${\tiny (0.4263,0.5754)}$ \\ \hline
& ${\tiny b=0.7}$ & ${\tiny 0.7015377}$ & ${\tiny 0.0015377}$ & ${\tiny %
0.0028327}$ & ${\tiny 0.0028351}$ & ${\tiny (0.6049,0.8093)}$ \\ \hline
& ${\tiny \theta }_{{\tiny 1}}{\tiny =0.8}$ & ${\tiny 0.8486051}$ & ${\tiny %
0.0486051}$ & ${\tiny 0.0765937}$ & ${\tiny 0.0789561}$ & ${\tiny %
(0.4596,1.5080)}$ \\ \hline
& ${\tiny \theta }_{{\tiny 2}}{\tiny =1.2}$ & ${\tiny 1.2433257}$ & ${\tiny %
0.0433257}$ & ${\tiny 0.1115477}$ & ${\tiny 0.1134248}$ & ${\tiny %
(0.7225,1.9904)}$ \\ \hline
& ${\tiny \theta }_{{\tiny 3}}{\tiny =1.3}$ & ${\tiny 1.3260033}$ & ${\tiny %
0.0260033}$ & ${\tiny 0.0700138}$ & ${\tiny 0.0706899}$ & ${\tiny %
(0.8795,1.9229)}$ \\ \hline\hline
200 & $\alpha =1.5$ & ${\tiny 1.5029588}$ & ${\tiny 0.0029588}$ & ${\tiny %
0.0057007}$ & ${\tiny 0.0057094}$ & ${\tiny (1.3728,1.6699)}$ \\ \hline
& ${\tiny a=0.5}$ & ${\tiny 0.4999514}$ & ${\tiny -.0000486}$ & ${\tiny %
0.0007961}$ & ${\tiny 0.0007961}$ & ${\tiny (0.4468,0.5564)}$ \\ \hline
& ${\tiny b=0.7}$ & ${\tiny 0.700456}$ & ${\tiny 0.000456}$ & ${\tiny %
0.0015143}$ & ${\tiny 0.0015145}$ & ${\tiny (0.6319,0.7807)}$ \\ \hline
& ${\tiny \theta }_{{\tiny 1}}{\tiny =0.8}$ & ${\tiny 0.8273891}$ & ${\tiny %
0.0273891}$ & ${\tiny 0.0397173}$ & ${\tiny 0.0404674}$ & ${\tiny %
(0.5076,1.2493)}$ \\ \hline
& ${\tiny \theta }_{{\tiny 2}}{\tiny =1.2}$ & ${\tiny 1.2235063}$ & ${\tiny %
0.0235063}$ & ${\tiny 0.0584567}$ & ${\tiny 0.0590093}$ & ${\tiny %
(0.8101,1.7279)}$ \\ \hline
& ${\tiny \theta }_{{\tiny 3}}{\tiny =1.3}$ & ${\tiny 1.3136215}$ & ${\tiny %
0.0136215}$ & ${\tiny 0.0367075}$ & ${\tiny 0.036893}$ & ${\tiny %
(0.9755,1.7035)}$ \\ \hline\hline
\end{tabular}%
\end{equation*}%
\begin{equation*}
\begin{tabular}{|c|c|c|c|c|c|c|}
\multicolumn{7}{l}{$\text{\textbf{Table 10.} The MLEs, Bias, Var, MSE and C.I%
\textbf{\ }values}$} \\ \hline\hline
n & parameter & Estimate & Bias & Var & MSE & C.I \\ \hline\hline
30 & ${\tiny \alpha =2}$ & ${\tiny 1.8540053}$ & ${\tiny -0.1459947}$ & $%
{\tiny 0.0861492}$ & ${\tiny 0.1074636}$ & ${\tiny (1.3579,2.5402)}$ \\ 
\hline
& ${\tiny a=0.2}$ & ${\tiny 0.2043051}$ & ${\tiny 0.0043051}$ & ${\tiny %
0.0047227}$ & ${\tiny 0.0047412}$ & ${\tiny (0.0793,0.3422)}$ \\ \hline
& ${\tiny b=1.5}$ & ${\tiny 1.5124662}$ & ${\tiny 0.0124662}$ & ${\tiny %
0.0168988}$ & ${\tiny 0.0170542}$ & ${\tiny (1.2850,1.7882)}$ \\ \hline
& ${\tiny \theta }_{1}{\tiny =0.5}$ & ${\tiny 0.5857591}$ & ${\tiny 0.0857591%
}$ & ${\tiny 0.092056}$ & ${\tiny 0.0994106}$ & ${\tiny (0.2091,1.3758)}$ \\ 
\hline
& ${\tiny \theta }_{2}{\tiny =0.6}$ & ${\tiny 0.6827111}$ & ${\tiny 0.0827111%
}$ & ${\tiny 0.108878}$ & ${\tiny 0.1157192}$ & ${\tiny (0.2713,1.4742)}$ \\ 
\hline
& ${\tiny \theta }_{3}{\tiny =0.9}$ & ${\tiny 0.9424471}$ & ${\tiny 0.0424471%
}$ & ${\tiny 0.0855705}$ & ${\tiny 0.0873723}$ & ${\tiny (0.4933,1.6119)}$
\\ \hline\hline
50 & ${\tiny \alpha =2}$ & ${\tiny 1.8900924}$ & ${\tiny -0.1099076}$ & $%
{\tiny 0.054245}$ & ${\tiny 0.0663247}$ & ${\tiny (1.4403,2.4056)}$ \\ \hline
& ${\tiny a=0.2}$ & ${\tiny 0.2016888}$ & ${\tiny 0.0016888}$ & ${\tiny %
0.0029133}$ & ${\tiny 0.0029162}$ & ${\tiny (0.0983,0.3127)}$ \\ \hline
& ${\tiny b=1.5}$ & ${\tiny 1.5054807}$ & ${\tiny 0.0054807}$ & ${\tiny %
0.0103271}$ & ${\tiny 0.0103572}$ & ${\tiny (1.3217,1.7219)}$ \\ \hline
& ${\tiny \theta }_{1}{\tiny =0.5}$ & ${\tiny 0.5604422}$ & ${\tiny 0.0604422%
}$ & ${\tiny 0.0556656}$ & ${\tiny 0.0593188}$ & ${\tiny (0.2489,1.1597)}$
\\ \hline
& ${\tiny \theta }_{2}{\tiny =0.6}$ & ${\tiny 0.657993}$ & ${\tiny 0.057993}$
& ${\tiny 0.0633326}$ & ${\tiny 0.0666958}$ & ${\tiny (0.3180,1.3028)}$ \\ 
\hline
& ${\tiny \theta }_{3}{\tiny =0.9}$ & ${\tiny 0.9271563}$ & ${\tiny 0.0271563%
}$ & ${\tiny 0.0536466}$ & ${\tiny 0.054384}$ & ${\tiny (0.5486,1.4589)}$ \\ 
\hline\hline
100 & ${\tiny \alpha =2}$ & ${\tiny 1.9448174}$ & ${\tiny -0.0551826}$ & $%
{\tiny 0.0271835}$ & ${\tiny 0.0302286}$ & ${\tiny (1.6074,2.2876)}$ \\ 
\hline
& ${\tiny a=0.2}$ & ${\tiny 0.2006911}$ & ${\tiny 0.0006911}$ & ${\tiny %
0.0014431}$ & ${\tiny 0.0014436}$ & ${\tiny (0.1273,0.2765)}$ \\ \hline
& ${\tiny b=1.5}$ & ${\tiny 1.5025237}$ & ${\tiny 0.0025237}$ & ${\tiny %
0.0050839}$ & ${\tiny 0.0050903}$ & ${\tiny (1.3686,1.6470)}$ \\ \hline
& ${\tiny \theta }_{1}{\tiny =0.5}$ & ${\tiny 0.5351945}$ & ${\tiny 0.0351945%
}$ & ${\tiny 0.0261092}$ & ${\tiny 0.0273479}$ & ${\tiny (0.2977,0.9242)}$
\\ \hline
& ${\tiny \theta }_{2}{\tiny =0.6}$ & ${\tiny 0.6309327}$ & ${\tiny 0.0309327%
}$ & ${\tiny 0.0311299}$ & ${\tiny 0.0320867}$ & ${\tiny (0.3745,1.0547)}$
\\ \hline
& ${\tiny \theta }_{3}{\tiny =0.9}$ & ${\tiny 0.9152326}$ & ${\tiny 0.0152326%
}$ & ${\tiny 0.0263485}$ & ${\tiny 0.0265805}$ & ${\tiny (0.6396,1.2798)}$
\\ \hline\hline
200 & ${\tiny \alpha =2}$ & ${\tiny 1.977254}$ & ${\tiny -0.022746}$ & $%
{\tiny 0.0112714}$ & ${\tiny 0.0117888}$ & ${\tiny (1.7652,2.1946)}$ \\ 
\hline
& ${\tiny a=0.2}$ & ${\tiny 0.2002004}$ & ${\tiny 0.0002004}$ & ${\tiny %
0.0006825}$ & ${\tiny 0.0006825}$ & ${\tiny (0.1481,0.2525)}$ \\ \hline
& ${\tiny b=1.5}$ & ${\tiny 1.5009306}$ & ${\tiny 0.0009306}$ & ${\tiny %
0.0023884}$ & ${\tiny 0.0023893}$ & ${\tiny (1.4047,1.6006)}$ \\ \hline
& ${\tiny \theta }_{1}{\tiny =0.5}$ & ${\tiny 0.5172642}$ & ${\tiny 0.0172642%
}$ & ${\tiny 0.0120296}$ & ${\tiny 0.0123276}$ & ${\tiny (0.3456,0.7748)}$
\\ \hline
& ${\tiny \theta }_{2}{\tiny =0.6}$ & ${\tiny 0.6161236}$ & ${\tiny 0.0161236%
}$ & ${\tiny 0.014622}$ & ${\tiny 0.014882}$ & ${\tiny (0.4252,0.9035)}$ \\ 
\hline
& ${\tiny \theta }_{3}{\tiny =0.9}$ & ${\tiny 0.9085774}$ & ${\tiny 0.0085774%
}$ & ${\tiny 0.0124857}$ & ${\tiny 0.0125593}$ & ${\tiny (0.7099,1.1482)}$
\\ \hline\hline
\end{tabular}%
\end{equation*}%
From Tables 9 and 10, we note that the bias is reduced as the sample size is
increased.

\section{Conclusions}

In this paper, we have proposed a bivariate exponentiated generalized linear
exponential distribution (BEGLED), whose marginals are exponentiated
generalized linear exponential distributions. We discussed some statistical
and reliability properties of the new distribution. Since the joint CDF and
the joint PDF are in a closed form, therefore the BEGLED can be used in
practice for non-negative and positively correlated random variables. The
maximum likelihood estimates (MLE) of the six parameters index to the BEGLED
are discussed. Moreover, a real data set is analyzed to show the usefulness
of the proposed distribution. Also, the bias of the parameters is calculated
using simulation studies. We hope our new distribution (BEGLED) might
attract wider sets of applications in reliability analysis.

\end{document}